\documentclass[letterpaper, 10 pt, conference]{ieeeconf}  

\IEEEoverridecommandlockouts                              
\usepackage{url}  
\usepackage{graphicx}  
\usepackage{amsmath,amssymb,amsfonts,mathrsfs,bm} 
\usepackage{caption}
\usepackage{cleveref}
\usepackage{array}
\usepackage{scalerel}
\usepackage{algorithm}
\usepackage{algpseudocode}
\usepackage{booktabs}
\usepackage{epstopdf}
\usepackage{multirow}
\usepackage{rotating}
\usepackage{algpseudocode}
\usepackage{cite}
\usepackage{epsfig}
\usepackage{tabularx}
\usepackage[table]{xcolor}
\usepackage{epstopdf}
\usepackage{accents}
\usepackage{circuitikz}
\usepackage{textcomp}
\usepackage{threeparttable,multirow}
\usetikzlibrary{arrows,shapes,positioning}
\usepackage{tikz}
\usepackage{pgfplots}
\definecolor{Fcolor}{rgb}{0, 0.5, 0.25}
\usepackage{subcaption}
\usepackage{caption}

\ifx\theorem\undefined
\newtheorem{theorem}{Theorem}
\newenvironment{theorem*}{\par\noindent{\bf Theorem\ }}{\hfill\\[2mm]}

\newenvironment{corollary*}{\par\noindent{\bf Corollary\ }}{\hfill\\[2mm]}
\newtheorem{definition}[theorem]{Definition}

\fi

\newcommand{\tr}{\mathrm{tr}}
\newcommand{\onebf}{\mathbf{1}}
\newcommand{\msh}{\!\!\:}

\newcommand{\Abf}{\mathbf{A}}

\newcommand{\bbf}{\mathbf{b}}

\newcommand{\Cbf}{\mathbf{C}}
\newcommand{\cbf}{\mathbf{c}}
\newcommand{\Ccal}{\mathcal{C}}
\newcommand{\Cbb}{\mathbb{C}}

\newcommand{\dbf}{\mathbf{d}}

\newcommand{\ebf}{\mathbf{e}}

\newcommand{\Ecal}{\mathcal{E}}

\newcommand{\fbf}{\mathbf{f}}

\newcommand{\Gcal}{\mathcal{G}}

\newcommand{\Hbf}{\mathbf{H}}

\newcommand{\Hbb}{\mathbb{H}}

\newcommand{\Ibf}{\mathbf{I}}
\newcommand{\irm}{\mathrm{i}}

\newcommand{\Mbf}{\mathbf{M}}

\newcommand{\Ncal}{\mathcal{N}}

\newcommand{\obf}{\mathbf{o}}

\newcommand{\pbf}{\mathbf{p}}

\newcommand{\qbf}{\mathbf{q}}

\newcommand{\Rbb}{\mathbb{R}}

\newcommand{\sbf}{\mathbf{s}}

\newcommand{\Tcal}{\mathcal{T}}

\newcommand{\ubf}{\mathbf{u}}

\newcommand{\Ucal}{\mathcal{U}}

\newcommand{\rbf}{\mathbf{r}}
\newcommand{\vbf}{\mathbf{v}}

\newcommand{\Vcal}{\mathcal{V}}

\newcommand{\Wbf}{\mathbf{W}}

\newcommand{\xbf}{\mathbf{x}}

\newcommand{\ybf}{\mathbf{y}}
\newcommand{\Ybf}{\mathbf{Y}}

\newcommand{\zbf}{\mathbf{z}}

\newcommand{\diagrm}{\mathrm{diag}}

\newcommand{\sm}[2]{\scaleto{#1\mathstrut}{#2pt}}
\newcommand{\smbullet}{\sm{\bullet}{5}}
\newif\ifcomment
\newcommand{\smallMinus}{\scalebox{1.2}[0.7]{-}}

\crefrangelabelformat{equation}{(#3#1#4)\,--\,(#5#2#6)}
\crefmultiformat{equation}{(#2#1#3)}{\,--\,(#2#1#3)}{#2#1#3}{#2#1#3}

\usepackage{xspace}
\newcommand{\matpower}[0]{\textsc{Matpower}\xspace}

\crefname{equation}{}{}
\crefname{figure}{Figure}{Figures}
\crefname{algorithm}{Algorithm}{}
\crefname{table}{Table}{Tables}
\crefname{lemma}{Lemma}{Lemmas}
\crefname{theorem}{Theorem}{Theorems}
\crefname{section}{Section}{Sections}
\crefname{definition}{Definition}{Definitions}
\commentfalse

\newcommand\ubar[1]{\underaccent{\bar}{#1}}
\makeatletter
\DeclareRobustCommand{\cev}[1]{%
  \mathpalette\do@cev{#1}%
}
\newcommand{\do@cev}[2]{%
  \fix@cev{#1}{+}%
  \reflectbox{$\m@th#1\vec{\reflectbox{$\fix@cev{#1}{-}\m@th#1#2\fix@cev{#1}{+}$}}$}%
  \fix@cev{#1}{-}%
}
\newcommand{\fix@cev}[2]{%
  \ifx#1\displaystyle
    \mkern#23mu
  \else
    \ifx#1\textstyle
      \mkern#23mu
    \else
      \ifx#1\scriptstyle
        \mkern#22mu
      \else
        \mkern#22mu
      \fi
    \fi
  \fi
}
\title{\LARGE \bf
Sequential Relaxation of Unit Commitment with\\ AC Transmission Constraints
}
\author{Fariba Zohrizadeh, Mohsen Kheirandishfard, Adnan Nasir, and Ramtin Madani
\thanks{Fariba Zohrizadeh and Mohsen Kheirandishfard are with the Department of Computer Science and Engineering, the University of Texas at Arlington, (email: fariba.zohrizadeh@uta.edu, mohsen.kheirandishfard@uta.edu). Adnan Nasir and Ramtin Madani are with the Department of Electrical Engineering, the University of Texas at Arlington (email: adnan.nasir@mavs.uta.edu, ramtin.madani@uta.edu). This work is in part supported by the NSF award 1809454 and a University of Texas System STARs award.
}%
}
\begin{document}

\maketitle
\thispagestyle{empty}
\pagestyle{empty}

\begin{abstract}
This paper proposes a sequential convex relaxation method for obtaining feasible and near-globally optimal solutions for unit commitment (UC) with AC transmission constraints. First, we develop a second-order cone programming (SOCP) relaxation for AC unit commitment. To ensure that the resulting solutions are feasible for the original non-convex problem, we incorporate penalty terms into the objective of the proposed SOCP relaxation. We generalize our penalization method to a sequential algorithm which starts from an initial point (not necessarily feasible) and leads to feasible and near-optimal solutions for AC unit commitment. Once a feasible point is attained, the algorithm preserves feasibility and improves the objective value until a near optimal point is obtained. The experimental results on IEEE 57, IEEE 118, and IEEE 300 bus benchmark cases from \matpower \cite{zimmerman2011matpower} demonstrate the performance of the proposed method in solving challenging instances of AC unit commitment.
\end{abstract}

\section{Introduction}
The unit commitment (UC) is a classical problem in the area of power systems which involves determining the optimal schedule for power generating units throughout a given planning horizon. The main objective is to meet power demand with minimum production cost while respecting the limitations of generating units and network constraints. Due to the economic importance of the UC problem, it has been heavily investigated for decades and is proven to be computationally hard in general \cite{tseng1996power,guan2003optimization}. 
The reader is referred to \cite{saravanan2013solution,allen2012price} and the references therein, for detailed surveys of the conventional formulations and methods for solving unit commitment.

A general unit commitment problem can be formulated as a mixed-integer optimization whose solution specifies the optimal status of generating units as well as voltages and power flows throughout the planning horizon. 
Additionally, several papers have considered uncertainties of demand and renewable generation into consideration using stochastic and robust optimization frameworks \cite{bitar2012bringing,bertsimas2013adaptive,phan2013optimization,yu2015impacts,lorca2017multistage,zhao2017unit,sundar2017modified,zheng2015stochastic}.
The incorporation of several other power system optimization problems into unit commitment has been envisioned as well,
such as the optimal power flow \cite{bai2009semi,castillo2016unit,lipka2017running}, network topology control \cite{hedman2010co}, demand response \cite{wu2013hourly}, air quality control \cite{kerl2015new}, and scheduling of deferrable loads \cite{subramanian2013real}. 

Various optimization methods have been used to approach the UC problem, such as branch-and-bound techniques \cite{cohen1983branch,marcovecchio2014deterministic,knueven2017ramping,dillon1978integer,carrion2006computationally,delarue2008adaptive,muckstadt1968application} and convex relaxations \cite{li2005price,knueven2017novel,finardi2014comparative}. 
In order to improve the efficiency of branch-and-bound searches, many papers have offered partial convex hull characterizations of UC feasible sets \cite{ostrowski2012tight,geng2017alternative,damci2016polyhedral,lee2004min}. Conic inequalities are proposed in \cite{akturk2009strong,frangioni2009computational,jabr2013rank,bai2009semi} to strengthen convex relaxations in the presence of nonlinear cost functions. 
In \cite{paredes2015using}, a combination of semidefinite programming relaxation and branch-and-bound is used to solve the day-ahead hydro unit commitment problem. In \cite{fattahi2017conic,ashraphijuo2016strong}, reformulation-linearization cuts are proposed to strengthen semidefinite programming relaxations of unit commitment. 
In \cite{bai2015decomposition}, a decomposition method is developed based on second-order cone programming (SOCP) to solve network constrained unit commitment with AC power flow constraints. In \cite{quan2014tighter}, a family of valid inequalities are proposed to improve the quality of SOCP relaxations of unit commitment. 
In \cite{liu2018global}, a global search algorithm is proposed which solves a sequence of mixed-integer second-order cone programming (MISOCP) problems, as well as nonlinear non-convex problems to lower- and upper-bound the globally optimal cost of unit commitment. 
In \cite{kargarian2015distributed,papavasiliou2015applying}, distributed frameworks on high-performance computing platforms are investigated for solving large-scale UC problems. Nevertheless, the improvements in run-time are reported to diminish with more than 15 parallel workers \cite{papavasiliou2013comparative}. 

In this paper, we introduce a novel sequential convex relaxation for solving unit commitment with AC transmission constraints. We propose a penalization method which is guaranteed to recover feasible solutions for general non-convex optimization problems under certain assumptions \cite{BMI1,QCQP_conic}. The proposed penalized convex relaxation can be solved sequentially in order to find feasible and near-globally optimal solutions. Our experimental results verify the effectiveness of this procedure in solving AC unit commitment problems on IEEE 57, IEEE 118, and IEEE 300 bus benchmark systems.
 
\subsection{Notations}
Throughout this paper, matrices, vectors, and scalars are represented by boldface uppercase, boldface lowercase, and italic lowercase letters, respectively. The symbols $\Rbb$, $\Cbb$, and $\Hbb_n$ denote the sets of real numbers, complex numbers, and $n\times n$ Hermitian matrices, respectively. The notation ``$\irm$'' is reserved for the imaginary unit. Notation $|\cdot|$ denotes either the absolute value of a scalar or the cardinality of a set, depending on the context. The symbols $(\cdot)^\ast$ and $(\cdot)^\top$ represent the conjugate transpose and transpose operators, respectively. For a given matrix $\Abf$, the notations $\Abf_{\smbullet,k}$, $\Abf_{j,\smbullet}$, and $A_{jk}$ refer to the $k^{th}$ column, $j^{th}$ row, and $(j,k)^{th}$ entry of the matrix $\Abf$, respectively. The Notation $\Abf\succeq 0$ means that $\Abf$ is symmetric/Hermitian and positive semidefinite.
\begin{table}[t]
\vspace{0.3cm}
\line(1,0){245}\\
\vspace{-0.4cm}
\begin{flushleft}
	\hspace{0.2cm}{\bf Unit Constraints:}
\end{flushleft}
\vspace{-0.6cm}

\begin{subequations}
\begin{align}
& x_{g,t}\in\{0,1\}\label{eq:UC_cons1}\\[1.5pt]
&c_{g,t}\!=\!{\alpha_{g}p_{g,t} \!+\!\beta_{g} p_{g,t}^{2}}+\notag\\& \phantom{c_{g,t}\!=\!}{{\gamma}_{g}x_{g,t}}+\!\gamma_{g}^{\uparrow}(1\!-\! x_{g,t\smallMinus 1})x_{g,t}\! + \!\gamma_{g}^{\downarrow} x_{g,t\smallMinus 1}(1 \!-\! x_{g,t} )\label{eq:UC_cons2}\\
& x_{g,\tau}\!-\! x_{g,\tau\smallMinus 1}\! \leq \! x_{g,t} &&&&\hspace{-3.5cm} \forall {\tau\!\in\!\{t\sm{-}{9} m_{g}^{\sm{\uparrow}{5}}\sm{+}{9}1, \dots, t\}}\label{eq:UC_cons3}\\[1.5pt]
& x_{g,\tau\smallMinus 1}\!-\! x_{g,\tau}\! \leq \!1\!-\! x_{g,t}&&&&\hspace{-3.5cm} \forall {\tau\!\in\!\{t\sm{-}{9} m_{g}^{\sm{\downarrow}{5}}\sm{+}{9}1, \dots, t\}}\label{eq:UC_cons4}\\[1.5pt]
& {\ubar{p}_{g}}\,{x_{g,t}} \leq {p_{g,t}} \leq {\bar{p}_{g}}\,{x_{g,t}}\label{eq:UC_cons5}\\[1.5pt]
& {\ubar{q}_{g}}\,{x_{g,t}} \leq {q_{g,t}} \!\: \leq {\bar{q}_{g}}\,{x_{g,t}} \label{eq:UC_cons6}\\[1.5pt]
& p_{g,t}\!-\!{p_{g,t\smallMinus 1}}\, \leq {r_{g}}{x_{g,t\smallMinus 1}}\!+\!{s_{g}}(1-x_{g,t\smallMinus 1})\label{eq:UC_cons7}\\[1.5pt]
& {p_{g,t\smallMinus1}}\!-\!{p_{g,t}}\, \leq {r_{g}}{x_{g,t}}\!+\!{s_{g}}(1-{x_{g,t}})\label{eq:UC_cons8}
\end{align}
\end{subequations}
\line(1,0){245}\\
\vspace{-0.4cm}
\begin{flushleft}
	\hspace{0.2cm}{\bf AC Network Constraints:}
\end{flushleft}
\vspace{-0.6cm}

\begin{subequations}
\begin{align}
&\hspace{-0.25cm}\dbf_{\smbullet,t}\!+\!\mathrm{diag}\{\vbf_{\smbullet,t}^{\phantom{\ast}}\vbf_{\smbullet,t}^{\ast}\Ybf^{\ast}\}\,=\Cbf^{\top}\!(\pbf_{\smbullet,t}+\!\irm\!\:\qbf_{\smbullet,t})\label{eq:NC_cons1}&&&&\\[1.5pt]
&\hspace{-0.25cm}\mathrm{diag}\{\vec{\Cbf}^{\phantom{\ast}}\vbf_{\smbullet,t}^{\phantom{\ast}}\vbf_{\smbullet,t}^{\ast}\vec{\Ybf}^{\ast}\} =\vec{\sbf}_{\smbullet,t}\label{eq:NC_cons2}&&\\[1.5pt]
&\hspace{-0.25cm}\mathrm{diag}\{\cev{\Cbf}^{\phantom{\ast}}{\vbf_{\smbullet,t}^{\phantom{\ast}}\vbf_{\smbullet,t}^{\ast}}\cev{\Ybf}^{\ast}\}= \cev{\sbf}_{\smbullet,t}\label{eq:NC_cons3}&&\\[1.5pt]
&\hspace{-0.25cm}\ubar{\vbf}\leq\lvert\vbf_{\smbullet,t}\rvert\leq \bar{\vbf}\label{eq:NC_cons4}\\[1.5pt]
&\hspace{-0.25cm}|\vec{\sbf}_{\smbullet,t}|^{2} \leq \fbf_{\mathrm{max};t}^2\label{eq:NC_cons5}\\[1.5pt]
&\hspace{-0.25cm}|\cev{\sbf}_{\smbullet,t}|^{2} \leq \fbf_{\mathrm{max};t}^2\label{eq:NC_cons6}
\end{align}
\end{subequations}
\vspace{-1mm}
\line(1,0){245}
	\caption{Unit and network constraints in power system scheduling.}
\label{fig:table_cons}
	
\vspace{-0mm}
\end{table}
\section{Problem Formulation}
The unit commitment (UC) problem aims at finding the most reliable and cost-efficient schedule for a set of generating units throughout a discrete time horizon $\Tcal$, subject to forecasted electricity demands and operational constraints. Let $\Gcal$ denote the set of generating units whose schedule needs to be determined. Define $x_{g,t}\in\{0,1\}$ as a binary variable indicating whether the generating unit $g\in\Gcal$ is committed during the time slot $t\in\Tcal$. If $x_{g,t}=1$, the unit is active and generates power within its capacity limitations, otherwise, no power is produced by $g$ during the time interval $t$. Define $p_{g,t}$ and $q_{g,t}$, respectively, as the amounts of active power and reactive power injections of generator $g$ during the time interval $t$. 

Denoted $\Vcal$ and $\Ecal$ as the sets of buses and branches in the network, respectively. For every bus $k\in\Vcal$, the demand forecast at time $t$ is denoted as $d_{k,t}\in\Cbb$, whose real and imaginary parts account for active and reactive power demands, respectively. Let $\Cbf\in\{0,1\}^{|\Gcal|\times|\Vcal|}$ be the incidence matrix whose $(g,k)$ entry is equal to $1$, if and only if the generating unit $g$ belongs to the bus $k$. Define the matrices $\vec{\Cbf}, \cev{\Cbf}\in\{0,1\}^{|\Ecal|\times|\Vcal|}$ as the \textit{from} and \textit{to} incidence matrices, respectively. The $(l,k)$ entry of $\vec{\Cbf}$ is equal to one, if and only if the line $l\in\Ecal$ starts at bus $k$, while the $(l,k)$ entry of $\cev{\Cbf}$ is equal to $1$, if and only if the line $l$ ends at bus $k$. Additionally, define $\Ybf\in\Cbb^{|\Vcal|\times|\Vcal|}$ as the nodal admittance matrices of the network and $\vec{\Ybf}$, ${\cev{\Ybf}}\in\Cbb^{|\Ecal|\times|\Vcal|}$ as the \textit{from} and \textit{to} branch admittance matrices. 

The feasible set of AC unit commitment can be described by unit constraints and AC network constraints. Unit constraints impose the minimum up and down time limits \cref{eq:UC_cons3,eq:UC_cons4}, generator capacities \cref{eq:UC_cons5,eq:UC_cons6}, as well as ramp limits \cref{eq:UC_cons7,eq:UC_cons8}. Define $m^{\uparrow}_{g}$ and $m^{\downarrow}_{g}$, respectively, as the minimum up time and minimum down time limits for generating unit $g$. 
If the unit $g$ is committed during the interval $t$, then its, active and reactive power injections must lie within the intervals $[\ubar{p}_{g},\bar{p}_{g}]$ and $[\ubar{q}_{g},\bar{q}_{g}]$, respectively. Additionally, denote $r_{g}$ as the maximum variation of active power injection by unit $g$ between two consecutive time slots in which the unit stays committed. Define $s_{g}$ as the maximum amount of active power injection after start-up and prior to shutdown. 

The network constraint \eqref{eq:NC_cons1} accounts for nodal power balances.
The constraint \eqref{eq:NC_cons4} enforces voltage magnitude limits.
Moreover, denote the line power flows at the starting and ending buses by $\vec\sbf\in\Cbb^{|\Ecal|\times |\Tcal|}$ and $\cev\sbf\in\Cbb^{|\Ecal|\times |\Tcal|}$, respectively. The constraints \eqref{eq:NC_cons5} and \eqref{eq:NC_cons6} enforce the thermal limits of lines. 

Given the above definitions, the AC unit commitment problem can be formulated as the optimization
\begin{subequations}
\begin{align}
&
{\text{minimize}}
&&\hspace{-0.3cm}\sum_{g,t} c_{g,t}\label{eq:UC_obj}\\
&\hspace{-0.0cm} \text{subject to}	&&\hspace{-0.3cm}(\xbf_{g,\smbullet}^{\!\top},\pbf_{g,\smbullet}^{\!\top},\qbf_{g,\smbullet}^{\!\top}, \cbf_{g,\smbullet}^{\!\top})\in\Ucal_{g}&&\forall g\in\Gcal,\!\!\label{eq:UC_set1}\\
&&&\hspace{-0.3cm}(\pbf_{\smbullet,t},\qbf_{\smbullet,t},\!\vbf_{\smbullet,t},\!\vec{\sbf}_{\smbullet,t},\cev{\sbf}_{\smbullet,t})\!\in\!\Ncal_{t}&&\forall t\in\Tcal\!,\label{eq:UC_set2}
\end{align}
\end{subequations}
with respect to the matrix variables $\xbf\triangleq [x_{g,t}]$, $\pbf \triangleq [p_{g,t}]$, $\qbf\triangleq[q_{g,t}]$, $\cbf\triangleq[c_{g,t}]$, $\vbf\triangleq[v_{k,t}]$, $\vec{\sbf}\triangleq[\vec{s}_{l,t}]$, and $\cev{\sbf}\triangleq[\cev{s}_{l,t}]$. The objective function \eqref{eq:UC_obj} is equal the sum of the production costs of all generating units throughout the time horizon $\Tcal$. For any arbitrary generating unit $g$ in time interval $t$, the production cost consists of the generation cost, start-up cost, shutdown cost, and a fixed cost. The generation cost is a quadratic function with respect to $p_{g,t}$ with nonnegative coefficients $\alpha_{g}$ and $\beta_{g}$. The start-up cost $\gamma^{\uparrow}_{g}$ and shutdown cost $\gamma^{\downarrow}_{g}$ are associated with every time slots at which the unit changes status. The fixed production cost $\gamma_{g}$ is enforced if the unit is active.  
\begin{definition}
\vspace{1mm}
For every generating units $g\in\Gcal$, define $\Ucal_{g}\!\subset\!\Rbb^{|\Tcal|\times 4}$ to be the set of all quadruplets $(\xbf_{g,\smbullet}^{\!\top},\pbf_{g,\smbullet}^{\!\top},\qbf_{g,\smbullet}^{\!\top}, \cbf_{g,\smbullet}^{\!\top})$ that satisfy the constraints \cref{eq:UC_cons1,eq:UC_cons2,eq:UC_cons3,eq:UC_cons4,eq:UC_cons5,eq:UC_cons6,eq:UC_cons7,eq:UC_cons8} throughout the entire planning horizon.
\vspace{1mm}
\end{definition}
\begin{definition}
\vspace{1mm}
For every $t\in\Tcal$, define $\Ncal_{t}\!\subset\!\Rbb^{|\Gcal|\times 2}$ $\!\times\Cbb^{|\Vcal|}\!\times\!\Cbb^{|\Ecal|\times 2}$ to be the set of all quintuplet $(\pbf_{\smbullet,t},\qbf_{\smbullet,t},\vbf_{\smbullet,t},\vec{\sbf}_{\smbullet,t},\cev{\sbf}_{\smbullet,t})$ that satisfy the network constraints \cref{eq:NC_cons1,eq:NC_cons2,eq:NC_cons3,eq:NC_cons4,eq:NC_cons5,eq:NC_cons6}.
\vspace{1mm}
\end{definition}

Problem \cref{eq:UC_obj,eq:UC_set1,eq:UC_set2} is a mixed-integer nonlinear optimization, due to the presence of binary variables and nonlinearity of the network constraints. In what follows, we will develop a convex relaxation to tackle the non-convexity of this problem.

\begin{table}[t]
\vspace{0.3cm}
\line(1,0){245}\\
\vspace{-0.4cm}
\begin{flushleft}
	\hspace{0.2cm}{\bf Unit Constraints:}
\end{flushleft}
\vspace{-0.6cm}

\begin{subequations}
\begin{align}
&z_{g,t}=x_{g,t},\label{eq:UC_RX_cons1}\\
&c_{g,t}\!=\!{\alpha_{g}p_{g,t} \!+\!\beta_{g} o_{g,t}}\!\!+\!{{\gamma}_{g}x_{g,t}}\notag\\& \phantom{c_{g,t}=\!}+\!\gamma_{g}^{\uparrow}(x_{g,t}\!-\! u_{g,t})\! + \!\gamma_{g}^{\downarrow} (x_{g,t\smallMinus 1} \!-\! u_{g,t} ),\label{eq:UC_RX_cons2}\\
& x_{g,\tau}\!-\! x_{g,\tau\smallMinus 1}\! \leq \! x_{g,t} &&&&\hspace{-3.5cm} \forall {\tau\!\in\!\{t\sm{-}{9} m_{g}^{\sm{\uparrow}{5}}\sm{+}{9}1, \dots, t\}}\label{eq:UC_RX_cons3}\\[1.5pt]
& x_{g,\tau\smallMinus 1}\!-\! x_{g,\tau}\! \leq \!1\!-\! x_{g,t}&&&&\hspace{-3.5cm} \forall {\tau\!\in\!\{t\sm{-}{9} m_{g}^{\sm{\downarrow}{5}}\sm{+}{9}1, \dots, t\}}\label{eq:UC_RX_cons4}\\[1.5pt]
& {\ubar{p}_{g}}\,{x_{g,t}} \leq {p_{g,t}} \leq {\bar{p}_{g}}\,{x_{g,t}}\label{eq:UC_RX_cons5}\\[1.5pt]
& {\ubar{q}_{g}}\,{x_{g,t}} \leq {q_{g,t}} \!\: \leq {\bar{q}_{g}}\,{x_{g,t}} \label{eq:UC_RX_cons6}\\[1.5pt]
& p_{g,t}\!-\!{p_{g,t\smallMinus 1}}\, \leq {r_{g}}{x_{g,t\smallMinus 1}}\!+\!{s_{g}}(1-x_{g,t\smallMinus 1})\label{eq:UC_RX_cons7}\\[1.5pt]
& {p_{g,t\smallMinus1}}\!-\!{p_{g,t}}\, \leq {r_{g}}{x_{g,t}}\!+\!{s_{g}}(1-{x_{g,t}})\label{eq:UC_RX_cons8}\\
        %
&\!\!\begin{bmatrix} z_{g,t-1} & u_{g,t}\\ u_{g,t} & z_{g,t} \end{bmatrix}-
\begin{bmatrix} x_{g,t-1} \\ x_{g,t} \end{bmatrix}
\begin{bmatrix} x_{g,t-1} & \hspace{-0.25cm}x_{g,t} \end{bmatrix}\!\succeq{0},\label{eq:UC_RX_cons9}\\
&\!\!\begin{bmatrix} z_{g,t\phantom{-1}} & b_{g,t}\\ b_{g,t\phantom{-1}} & o_{g,t}\, \end{bmatrix}-
\begin{bmatrix} \phantom{_-} x_{g,t\phantom{1}} \\ \phantom{_-} p_{g,t\phantom{1}} \end{bmatrix}
\begin{bmatrix} x_{g,t\phantom{-1}} &\hspace{-0.25cm}p_{g,t} \end{bmatrix}\!\succeq {0}\label{eq:UC_RX_cons10}.
\end{align}
\end{subequations}
\line(1,0){245}\\
\vspace{-0.4cm}
\begin{flushleft}
	\hspace{0.2cm}{\bf AC Network Constraints:}
\end{flushleft}
\vspace{-0.6cm}

\begin{subequations}
\begin{align}
&\hspace{-0.8cm}\dbf_{\smbullet,t}\!+\!\diagrm\{\Wbf_{t}^{\phantom\ast}\Ybf^{\ast}\}\,=\Cbf^{\top}\!(\pbf_{\smbullet,t}+\!\irm\!\:\qbf_{\smbullet,t}),\label{eq:NC_RX_cons1}&&\\[1.2pt]
&\hspace{-0.8cm}\diagrm\{\vec{\Cbf}^{\phantom\ast}\Wbf_{t}^{\phantom\ast}\vec{\Ybf}^{\ast}\}=\vec{\sbf}_{\smbullet,t}^{\phantom\ast},\label{eq:NC_RX_cons2}&&\\[1.2pt]
&\hspace{-0.8cm}\diagrm\{\cev{\Cbf}^{\phantom\ast}{\Wbf_{t}^{\phantom\ast}}\cev{\Ybf}^{\ast}\}= \cev{\sbf}_{\smbullet,t}^{\phantom\ast},\label{eq:NC_RX_cons3}&&\\[1.2pt]
&\hspace{-0.8cm}\ubar{\vbf}^{2}\leq\diagrm\{\Wbf_{t}\}\leq \bar{\vbf}^{2},\label{eq:NC_RX_cons4}\\[1.2pt]
&\hspace{-0.8cm}|\vec{\sbf}_{\smbullet,t}|^{2} \leq \vec{\fbf}_{\smbullet,t} \leq \fbf_{\mathrm{max};t}^2,\label{eq:NC_RX_cons5}\\[1.2pt]
&\hspace{-0.8cm}|\cev{\sbf}_{\smbullet,t}|^{2} \leq \cev{\fbf}_{\smbullet,t} \leq \fbf_{\mathrm{max};t}^2\label{eq:NC_RX_cons6}\\
&\hspace{-0.8cm}\Wbf_{t}^{\phantom\ast}-\vbf^{\phantom\ast}_{\smbullet,t}\vbf^{\ast}_{\smbullet,t}\succeq_{\,\sm{\Ccal}{6}}{0}.\label{eq:NC_RX_cons7}
\end{align}
\end{subequations}
\vspace{-1mm}
\line(1,0){245}
	\caption{Relaxed unit and AC network constraints.}
\label{fig:table_cons}
	
\vspace{-0mm}
\end{table}
\section{Convex Relaxation of the UC Problem}
The non-convex sets $\{\Ucal_{g}\}_{g\in\Gcal}$ and $\{\Ncal_{t}\}_{t\in\Tcal}$, are the sources of computational complexity. In this paper, we introduce convex surrogates $\{\Ucal_{g}^{\;\!\mathrm{conv}}\}_{g\in\Gcal}$ and $\{\Ncal_{t}^{\;\!\mathrm{conv}}\}_{t\in\Tcal}$, which lead to a class of computationally-tractable relaxations of the problem \cref{eq:UC_obj,eq:UC_set1,eq:UC_set2}. To this end, define the auxiliary variables $\ubf,\obf,\rbf,\zbf,\bbf\!\in\!\Rbb^{|\Gcal|\times|\Tcal|}$, whose components account for monomials $x_{g,t-1}x_{g,t}$, $p_{g,t}^{2}$, $q_{g,t}^{2}$, $x_{g,t}^{2}$, and $x_{g,t}p_{g,t}$, respectively. Using the defined variables, non-convex constraints \cref{eq:UC_cons1,eq:UC_cons2} can be convexified as \cref{eq:UC_RX_cons1,eq:UC_RX_cons2}. In addition, to relax the non-convexity of AC network constraints, we define the auxiliary variables $\vec\fbf_{\smbullet,t}, \cev\fbf_{\smbullet,t}\in\Rbb^{|\Ecal|}$ and $\Wbf_{\!t}\in\Hbb_{|\Vcal|}$, accounting for $|\vec{s}_{\smbullet,t}|^{2}$, $|\cev{s}_{\smbullet,t}|^{2}$, and $\vbf_{\smbullet,t}^{\phantom{\ast}}\vbf_{\smbullet,t}^{\ast}$, respectively. Using the above auxiliary variables, the non-convex constraints \cref{eq:NC_cons1,eq:NC_cons2,eq:NC_cons3,eq:NC_cons4,eq:NC_cons5} can be relaxed as \cref{eq:NC_RX_cons1,eq:NC_RX_cons2,eq:NC_RX_cons3,eq:NC_RX_cons4,eq:NC_RX_cons5,eq:NC_RX_cons6}.

In order to capture the binary requirements of the commitment decisions and enforce the relationship between the auxiliary variables and the corresponding monomials, we strengthen the proposed convex relaxation via conic constraints \cref{eq:UC_RX_cons9,eq:UC_RX_cons10}, and \eqref{eq:NC_RX_cons7}, where $\Ccal$ in \eqref{eq:NC_RX_cons7} is a pointed convex cone. Next, we define the convex surrogates $\{\Ucal_{g}^{\,\mathrm{conv}}\}_{g\in\Gcal}$ and $\{\Ncal_{t}^{\,\mathrm{conv}}\}_{t\in\Tcal}$.

\begin{definition}
\vspace{1mm}
For every $g\!\in\!\Gcal$, define $\Ucal_{g}^{\,\mathrm{conv}}\!\!\subset\Rbb^{|\Tcal|\times 9}$ to be the set of all nonuplets $(\xbf_{g,\smbullet}^{\!\top},\pbf_{g,\smbullet}^{\!\top},\qbf_{g,\smbullet}^{\!\top},\cbf_{g,\smbullet}^{\!\top},\ubf_{g,\smbullet}^{\!\top},\obf_{g,\smbullet}^{\!\top},\rbf_{g,\smbullet}^{\!\top},$ $\zbf_{g,\smbullet}^{\!\top},\bbf_{g,\smbullet}^{\!\top})$ that satisfy the constraints \cref{eq:UC_RX_cons1,eq:UC_RX_cons2,eq:UC_RX_cons3,eq:UC_RX_cons4,eq:UC_RX_cons5,eq:UC_RX_cons6,eq:UC_RX_cons7,eq:UC_RX_cons8,eq:UC_RX_cons9,eq:UC_RX_cons10} throughout the entire planning horizon.
\vspace{1mm}
\end{definition}

\begin{definition}
\vspace{1mm}
For every $t\!\in\!\Tcal$, define $\Ncal_{t}^{\,\mathrm{conv}}\!\subset\!\Rbb^{|\Gcal|\times 2} \times\!$ $\Cbb^{|\Vcal|} \!\times\! \Cbb^{|\Ecal|\times 2} \!\times\! \Hbb_{|\Vcal|} \!\times\! \Rbb^{|\Ecal|\times 2}$ to be the set of all octuplets $(\pbf_{\smbullet,t},\qbf_{\smbullet,t},\vbf_{\smbullet,t},\vec{\sbf}_{\smbullet,t},\cev{\sbf}_{\smbullet,t},\!\Wbf_{t},\vec\fbf_{\smbullet,t},\cev\fbf_{\smbullet,t})$ that satisfy the constraints \cref{eq:NC_RX_cons1,eq:NC_RX_cons2,eq:NC_RX_cons3,eq:NC_RX_cons4,eq:NC_RX_cons5,eq:NC_RX_cons6,eq:NC_RX_cons7}.
\vspace{1mm}
\end{definition}


The second-order cone programming (SOCP) relaxation of network constraints can be derived by incorporating the following convex set into the constraint \eqref{eq:NC_RX_cons7}:
\begin{equation}\label{eq:SOCP_Set}\notag
\begin{aligned}
\hat\Ccal\triangleq\big\{ \Hbf\in\Hbb_{|\Vcal|} \;\big|\; H_{ii}\geq\! 0,\;H_{ii}H_{jj}\geq|H_{ij}|^2\!\!,\;\forall (i,\!j)\!\in\!\Ecal \big\}.
\end{aligned}
\end{equation}
The solution provided by the SOCP relaxation is a lower-bound for the globally optimal solution of AC unit commitment.
In general, solutions obtained from convex relaxations are not necessarily feasible for the original non-convex problem. To remedy this shortcoming, we propose a novel penalization method to obtain feasible points. In the next section, we describe the proposed penalization method in details.
\section{Penalization Method}
We incorporate a linear penalty term $\kappa(\{\Wbf_{\!t}\}_{t\in\Tcal},\zbf,\obf,\rbf,$ $\vec{\fbf},\cev{\fbf},\vbf,\xbf,\sbf,\vec{\sbf},\cev{\sbf})$ into the objective of the relaxed problem to enforce feasibility. Given an initial guess $\ybf^{\sm{0}{5}}\msh=\msh (\vbf^{\sm{0}{5}}\msh,\xbf^{\sm{0}{5}}\msh,\sbf^{\sm{0}{5}}\msh,\vec{\sbf}^{\:\sm{0}{5}}\msh,\msh\cev{\sbf}^{\sm{0}{5}})$ that is sufficiently close to the feasible set of the problem \cref{eq:UC_obj,eq:UC_set1,eq:UC_set2}, the following choice of penalty function guarantees the feasibility of the resulting solution under the assumptions in \cite{BMI1,QCQP_conic}: 
\begin{align}
\small
&\kappa_{\Mbf,\ybf_0}(\{\Wbf_{\!t}\}_{t\in\Tcal},\zbf,\obf,\rbf,\vec{\fbf},\cev{\fbf},\vbf,\xbf,\sbf,\vec{\sbf},\cev{\sbf})\triangleq\nonumber\\
&\hspace{0.2cm}\sum\mathop{}_{\mkern-5mu t}(
\tr\{\Wbf_{\!t}\Mbf\}\!-\!{\vbf_{\smbullet,t}^{\sm{0}{5}}}^{\hspace{-0.15cm}\ast}\,\Mbf\vbf_{\smbullet,t}\!-\!
\vbf_{\smbullet,t}^{\ast}\,\Mbf\vbf_{\smbullet,t}^{\sm{0}{5}}+
{\vbf_{\smbullet,t}^{\sm{0}{5}}}^{\hspace{-0.15cm}\ast}\,\Mbf\vbf_{\smbullet,t}^{\sm{0}{5}}+\nonumber\\[1pt]
&\hspace{0.2cm}\phantom{\sum\mathop{}_{\mkern-5mu t}(}
\zbf_{\smbullet,t}^{\!\top}\onebf-2\,\xbf_{\smbullet,t}^{\!\top}\,{\xbf_{\smbullet,t}^{\sm{0}{5}}} + {\xbf_{\smbullet,t}^{\sm{0}{5}}}^{\hspace{-0.25cm}\top}\,\xbf_{\smbullet,t}^{\sm{0}{5}}+\nonumber\\[1pt]
&\hspace{0.2cm}\phantom{\sum\mathop{}_{\mkern-5mu t}(}
\obf_{\smbullet,t}^{\!\top}\onebf-2\,\pbf_{\smbullet,t}^{\!\top}\,\pbf_{\smbullet,t}^{\sm{0}{5}} + {\pbf_{\smbullet,t}^{\sm{0}{5}}}^{\hspace{-0.25cm}\top}\,\pbf_{\smbullet,t}^{\sm{0}{5}}+\nonumber\\[1pt]
&\hspace{0.2cm}\phantom{\sum\mathop{}_{\mkern-5mu t}(}
\rbf_{\smbullet,t}^{\!\top}\onebf-2\,\qbf_{\smbullet,t}^{\!\top}\,\qbf_{\smbullet,t}^{\sm{0}{5}} + {\qbf_{\smbullet,t}^{\sm{0}{5}}}^{\hspace{-0.25cm}\top}\,\qbf_{\smbullet,t}^{\sm{0}{5}}+\nonumber\\[1pt]
&\hspace{0.2cm}\phantom{\sum\mathop{}_{\mkern-5mu t}(}
\vec{\fbf}_{\smbullet,t}^{\,\top} \onebf -
{\vec{\sbf}_{\smbullet,t}^{\,\sm{0}{5}}}^{\hspace{-0.15cm}\ast}\:\vec{\sbf}_{\smbullet,t}-
{\vec{\sbf}_{\smbullet,t}}^{\hspace{-0.25cm}\ast}\:\,{\vec{\sbf}_{\smbullet,t}^{\,\sm{0}{5}}}+
{\vec{\sbf}_{\smbullet,t}^{\,\sm{0}{5}}}^{\hspace{-0.15cm}\ast}\:{\vec{\sbf}^{\,\sm{0}{5}}_{\smbullet,t}}+\nonumber\\[1pt]
&\hspace{0.2cm}\phantom{\sum\mathop{}_{\mkern-5mu t}(}
\cev{\fbf}_{\smbullet,t}^{\,\top} \onebf -
{\cev{\sbf}_{\smbullet,t}^{\,\sm{0}{5}}}^{\hspace{-0.15cm}\ast}\:\cev{\sbf}_{\smbullet,t}-
{\cev{\sbf}_{\smbullet,t}}^{\hspace{-0.25cm}\ast}\:\,{\cev{\sbf}_{\smbullet,t}^{\,\sm{0}{5}}}+
{\cev{\sbf}_{\smbullet,t}^{\,\sm{0}{5}}}^{\hspace{-0.15cm}\ast}\:{\cev{\sbf}^{\,\sm{0}{5}}_{\smbullet,t}}),\label{penpen}
\end{align}
where $\Mbf\in\Hbb_{|\Vcal|}$ is a fixed penalty matrix. 

By augmenting the penalty term \eqref{penpen} into the objective function of the relaxed problem, the \textit{penalized convex relaxation} of AC unit commitment can be formulated as:

\vspace{-0.3cm}
\begin{subequations}
\begin{align}
&\!\!\!\!
{\text{min}}
&&\hspace{-0.2cm}g(\cbf)\!+\!\mu\kappa_{\Mbf,\ybf_0}(\{\msh\Wbf_{\!t}\!\}_{t\in\Tcal}\!,\zbf,\msh\obf,\msh\rbf\msh,\msh\vec{\fbf}\!,\cev{\fbf}\!,\!\vbf\!,\msh\xbf,\msh\pbf\!+\!\irm \qbf,\vec{\sbf},\msh\cev{\sbf})\!\label{eq:UC_P_obj}\\
&\!\!\!\!\text{s.t.}	
&&\hspace{-0.2cm}(\xbf_{g,\smbullet}^{\!\top},\msh\pbf_{g,\smbullet}^{\!\top},\msh\qbf_{g,\smbullet}^{\!\top},  \msh\cbf_{g,\smbullet}^{\!\top},\msh\zbf_{g,\smbullet}^{\!\top},\msh\obf_{g,\smbullet}^{\!\top},\msh\rbf_{g,\smbullet}^{\!\top})\!\!\!\;\in\!\!\;\Ucal_{g}^{\,\mathrm{conv}}&&&&\hspace{-1.450cm}\forall g\!\in\!\Gcal\!\!\;,\label{eq:UC_P_set1}\\
&&&\hspace{-0.2cm}(\pbf_{\smbullet,t},\msh\qbf_{\smbullet,t},\!\vbf_{\smbullet,t},\!\!\;\vec{\sbf}_{\smbullet,t},\msh\cev{\sbf}_{\smbullet,t},\!\Wbf_{\!t},\msh\vec{\fbf}_{\smbullet,t},\msh\cev{\fbf}_{\smbullet,t})\!\!\!\;\in\!\Ncal_{t}^{\,\mathrm{conv}}&&&&\hspace{-1.450cm}\forall t\!\in\!\!\!\;\Tcal\!\!\;,\label{eq:UC_P_set2}
\end{align}
\end{subequations}
with respect to decision variables $\xbf\triangleq [x_{g,t}]$, $\pbf \triangleq [p_{g,t}]$, $\qbf\triangleq[q_{g,t}]$, $\cbf\triangleq[c_{g,t}]$, $\zbf\triangleq[z_{g,t}]$, $\obf\triangleq[o_{g,t}]$, $\rbf\triangleq[r_{g,t}]$, $\vbf\triangleq[v_{k,t}]$, $\vec{\sbf}\triangleq[\vec{s}_{l,t}]$, $\cev{\sbf}\triangleq[\cev{s}_{l,t}]$, $\vec{\fbf}\triangleq[\vec{f}_{l,t}]$, $\cev{\fbf}\triangleq[\cev{f}_{l,t}]$, and $\{\Wbf_{\!t}\}_{t\in\Tcal}$. The nonnegative penalty parameter $\mu>0$ sets the trade off between the objective and the penalty functions. The penalized convex relaxation \cref{eq:UC_P_obj,eq:UC_P_set1,eq:UC_P_set2} is said to be \textit{tight} if it possesses a unique optimal solution $(\xbf,\pbf,\qbf,\cbf,\zbf,\obf,\rbf,\vbf,\vec{\sbf},\cev{\sbf},\vec{\fbf},\cev{\fbf},\{\Wbf_{\!t}\}_{t\in\Tcal})$ such that $x_{g,t}\in\{0,1\}$ and $\Wbf_{t}=\vbf^{\phantom\ast}_{\smbullet,t}\vbf^{\ast}_{\smbullet,t}$, for every $g\in\Gcal$ and $t\in\Tcal$. The tightness of the penalization guarantees the recovery of a feasible point for AC unit commitment \cref{eq:UC_obj,eq:UC_set1,eq:UC_set2}.
\subsection{Choice of Penalty Matrix}
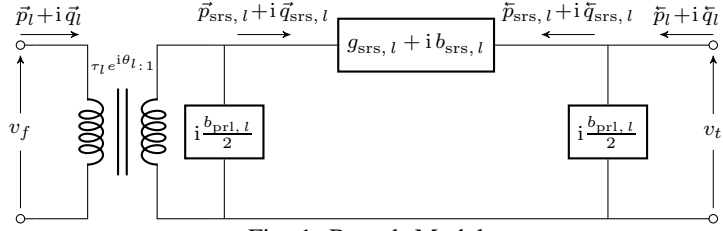
\begin{figure*}[!t]
  \begin{center}
  \vspace{2mm}
    \tikzset{R_F/.style = {draw,thick,minimum size=2em}}
    \begin{circuitikz}[american, cute inductors]
		\draw (0,0) node [scale = 1.1, transformer core] (T){}  (T.base) node[below] {\tiny$\tau_l e^{\irm\theta_l}\!\!:\!1$};
        \draw (T.A1) to [short, -o] ++(-0.2,0) coordinate (TA1);
	    \draw (T.A2) to [short, -o] ++(-0.2,0) coordinate (TA2);
	\path (TA2) -- node (A2A1) {\footnotesize $v_{f}$} (TA1);
		\draw[->,>=stealth',shorten >=4pt, shorten <=0.5pt, very thin] (A2A1) edge[bend left=0]  (TA1);
        \draw[-,shorten >=0.5pt, shorten <=4pt, very thin] (TA2) edge[bend left=0]  (A2A1);
    \draw[->,>=stealth',very thin] ([xshift=0cm,yshift=4pt]TA1) -- node[above] {\footnotesize $\vec{p}_{l}\!+\!\irm\,\vec{q}_{l}$} +(22pt,0);
    \draw (T.B1) to [short] ++(0.2,0) coordinate (B1E);
    \draw (T.B2) to [short] ++(0.2,0) coordinate (B2E);
    
    \path (B1E) -- node[R_F] (r1) {\footnotesize${
g_{\sm{\mathrm{srs}}{6},\,l} +\irm\,
b_{\sm{\mathrm{srs}}{6},\,l}}$} ++(5.1,0) coordinate (E1Y);
        	\draw[-] (r1) edge[bend left=0]  (E1Y);
        	\draw[-] (B1E) edge[bend left=0]  (r1);
            
    \draw (B2E) to [short] ++(5.1,0) coordinate (E2Y);
    \draw[->,>=stealth',very thin] ([xshift=5pt,yshift=4pt]B1E) -- node[above] {\footnotesize ${\vec{p}_{\sm{\mathrm{srs}}{6},\,l}\!+\! \irm \,\vec{q}_{\sm{\mathrm{srs}}{6},\,l}}$} +(20pt,0);
     \draw[->,>=stealth',very thin] ([xshift=-5pt,yshift=4pt]E1Y) -- node[above] {\footnotesize ${\cev{p}_{\sm{\mathrm{srs}}{6},\,l}\!+\! \irm \, \cev{q}_{\sm{\mathrm{srs}}{6},\,l}}$} +(-20pt,0);
    \draw (E1Y) to [short,-o] ++(1.4,0) coordinate (Y1E);
    \draw (E2Y) to [short,-o] ++(1.4,0) coordinate (Y2E);
    
	\path (B1E) -- node[R_F] (r2) {\footnotesize$\irm\frac{b_{\mathrm{prl},\,l}}{2}$} (B2E);
        	\draw[-] (r2) edge[bend left=0]  (B2E);
        	\draw[-] (B1E) edge[bend left=0]  (r2);
            
	\path (E1Y) -- node[R_F] (r3) {\footnotesize$\irm\frac{b_{\mathrm{prl},\,l}}{2}$} (E2Y);
        	\draw[-] (r3) edge[bend left=0]  (E2Y);
        	\draw[-] (E1Y) edge[bend left=0]  (r3);
            
    \draw[->,>=stealth',very thin] ([xshift=0cm,yshift=4pt]Y1E) -- node[above] {\footnotesize $\cev{p}_l\!+\! \irm \,\cev{q}_l$} +(-20pt,0);
	\path (Y2E) -- node (EY1Y2) {\footnotesize $v_{t}$} (Y1E);
        	\draw[->,>=stealth',shorten >=4pt, shorten <=1pt, very thin] (EY1Y2) edge[bend left=0]  (Y1E);
        	\draw[-,shorten >=1pt, shorten <=4pt, very thin] (Y2E) edge[bend left=0]  (EY1Y2);

   \end{circuitikz} 
  \end{center}
    \vspace{-4mm}
  \caption{Branch Model}
  \vspace{-4mm}
  \label{fig:branch}
\end{figure*}
Motivated by the previous literatures \cite{zohrizadeh2018opf,madani2016promises,madani2015convex}, we choose $\Mbf$ such that the penalty term $\tr\{\Wbf_{\!t}\Mbf\}$ reduces the apparent power loss over the series admittance of every line in the network. Consider the standard $\pi$-model of line $l\in\mathcal{E}$, with series admittance $y_{\sm{\mathrm{srs},\,l}{6}}\triangleq g_{\sm{\mathrm{srs},\,l}{6}}+\irm\,{b_{\sm{\mathrm{srs},\,l}{6}}}$ and total shunt susceptance $b_{\sm{\mathrm{prl},\,l}{6}}$, in series with a phase shifting transformer whose tap ratio has magnitude $\tau_l$ and phase shift angle $\theta_l$ \cite{zimmerman2011matpower}. 
The model is shown in Figure \ref{fig:branch}.
In order to penalize the apparent power loss over all lines of the network, we choose matrix $\Mbf$ as, 
\begin{equation}\notag
\begin{aligned}
\Mbf =\!\!\!\sum_{(i,j)\in\Ecal} [\ebf_{i}, \ebf_{j}](\Mbf_{ij}\!+\!\alpha\Ibf_{2})[\ebf_{i}, \ebf_{j}]^{\!\top}\!,
\end{aligned}
\end{equation}
where $\ebf_1,\ldots,\ebf_{|\Vcal|}$ denote the standard basis for $\Rbb^{|\Vcal|}$, and $\alpha$ is a positive constant. Moreover, each $\Mbf_{ij}$ is a $2\times 2$ positive semidefinite matrix defined as,
\begin{equation}\notag
\begin{aligned}
\Mbf_{ij} = 
\zeta_{ij}(\vec{\Ybf}_{\mathrm{q};\,l}+
\cev{\Ybf}_{\mathrm{q};\,l})
+\frac{\eta}{1\!-\!\eta}(\vec{\Ybf}_{\mathrm{p};\,l}+\cev{\Ybf}_{\mathrm{p};\,l}).
\end{aligned}
\end{equation}
where $\eta>0$ sets the trade-off between active and reactive loss minimization, and
%
\begin{equation}\notag
\begin{aligned}
&\!\vec{\Ybf}_{\mathrm{p};\,l}\!\triangleq\!\!\begin{bmatrix} \frac{g_{\sm{\mathrm{srs}}{6},\,l}}{\tau_l^2} & 
\!\!\frac{e^{\irm\theta_l}\,y_{\sm{\mathrm{srs}}{6},\,l}}{\smallMinus 2\tau_l}\! \\[1.0ex] 
\frac{\!y_{\sm{\mathrm{srs}}{6},\,l}^{\ast}}{\smallMinus2\tau_l e^{\irm\theta_l} } & 0 \\\vspace{-2.5mm}\end{bmatrix}\!\!,
&&\vec{\Ybf}_{\mathrm{q};\,l}\!\triangleq\!\!\begin{bmatrix} \frac{b_{\sm{\mathrm{srs}}{6},\,l}}{\smallMinus\tau_l^2} & \!\!\frac{ e^{\irm\theta_l}\,y_{\sm{\mathrm{srs}}{6},\,l}}{2\tau_l\irm} \\[1.0ex] \frac{ y_{\sm{\mathrm{srs}}{6},\,l}^{\ast}}{\smallMinus 2\tau_l\irm e^{\irm\theta_l}} & 0
\\\vspace{-2.5mm}\end{bmatrix}\!\!,\\
&\!\cev{\Ybf}_{\mathrm{p};\,l}\!\triangleq\!\!\begin{bmatrix} 0 &  
\frac{e^{\irm\theta_l}\,y_{\sm{\mathrm{srs}}{6},\,l}^{\ast}}{\smallMinus 2\tau_l} \\[1.0ex] \frac{ y_{\sm{\mathrm{srs}}{6},\,l}}{\smallMinus2\tau_l e^{\irm\theta_l}}\!\!\! & g_{\sm{\mathrm{srs}}{6},\,l}\\\vspace{-2.5mm}\end{bmatrix}\!\!, 
&&\cev{\Ybf}_{\mathrm{q};\,l}\!\triangleq\!\!\begin{bmatrix} 0 & \frac{ e^{\irm\theta_l}\,y_{\sm{\mathrm{srs}}{6},\,l}^{\ast}}{-2\tau_l\irm} \\[1.0ex]
\frac{ y_{\sm{\mathrm{srs}}{6},\,l}}{2\tau_l\irm e^{\irm\theta_l}} & \smallMinus b_{\sm{\mathrm{srs}}{6},\,l}\\\vspace{-2.5mm}\end{bmatrix}\!\!.
\end{aligned}
\end{equation}
Each $\zeta_{ij}\in\{-1,+1\}$ is determined based on the inductive or capacitive behavior of the line $l\!\in\!\Ecal$. More precisely, we set $\zeta_{ij}=1$ if the series admittance $y_{\sm{\mathrm{srs},\,l}{6}}$ is inductive (i.e., $b_{\sm{\mathrm{srs},\,l}{6}}\leq 0$), and $\zeta_{ij}=-1$, otherwise.

\subsection{Sequential Penalized Relaxation} 
The penalized SOCP relaxation \cref{eq:UC_P_obj,eq:UC_P_set1,eq:UC_P_set2} is guaranteed to produce a feasible solution for AC unit commitment if the initial guess $\ybf^{\sm{0}{5}}$ is sufficiently close to the feasible set of the original problem \cref{eq:UC_obj,eq:UC_set1,eq:UC_set2}. 
If a high quality initial point is not available, the proposed penalized SOCP relaxation can be solved sequentially 
until a feasible point for problem \cref{eq:UC_obj,eq:UC_set1,eq:UC_set2} is obtained. 
Once feasibility is attained, the sequential procedure 
improves the objective function while preserving the feasibility at each round until a near-optimal point is achieved. 
This sequential procedure is detailed by
\cref{al:alg_1}.
\begin{algorithm}[t]
\small
\caption{Sequential Penalized SOCP Relaxation.}\label{alg:alg_1}
\begin{algorithmic}[1]
\Require {$\mu$, $\Mbf$, $(\vbf_{\sm{0}{5}},\xbf_{\sm{0}{5}},\sbf_{\sm{0}{5}},\vec\sbf_{\sm{0}{5}},\cev\sbf_{\sm{0}{5}})$}
\Repeat
		\State Solve problem \cref{eq:UC_P_obj,eq:UC_P_set1,eq:UC_P_set2} to obtain $(\vbf,\xbf,\sbf,\vec\sbf,\cev\sbf)$
		\State $(\vbf_{\sm{0}{5}},\xbf_{\sm{0}{5}},\sbf_{\sm{0}{5}},\vec\sbf_{\sm{0}{5}},\cev\sbf_{\sm{0}{5}})\gets(\vbf,\xbf,\sbf,\vec\sbf,\cev\sbf)$
\Until {stopping criteria satisfied}
\Ensure {best found solution $(\vbf,\xbf,\sbf,\vec{\sbf},\cev{\sbf})$}
\end{algorithmic}\label{al:alg_1}
\end{algorithm}
\vspace{-0.1cm}
\section{Experimental Results} 
In this section, we present the results of our experiment on IEEE 57 bus, IEEE 118 bus, and IEEE 300 bus systems from \matpower \cite{zimmerman2011matpower}. The numerical experiments are performed in MATLAB using a 64-bit computer with an Intel 3.0 GHz, 12-core CPU, and 256 GB RAM. Note that the experiments are all performed on a workstation with a single CPU. The CVX package version 3.0 \cite{cvx} and MOSEK version 8.0 \cite{mosek2015mosek} are used to solve the proposed convex relaxations. 

The details of data generation are taken from \cite{madani2017scalable}. For each experiment, the cost coefficients $\alpha_g$, $\beta_g$, $\gamma_g$, $\gamma^{\downarrow}_g$ and $\gamma^{\uparrow}_g$ are chosen uniformly between zero and $\$1/(\mathrm{MW.h})^2$, $\$10/(\mathrm{MW.h})$, $\$100$, $\$30$ and $\$50$, respectively. The ramp limits of each generating unit are set to $r_g = s_g = \max\{\bar{p}_g/4,\ubar{p}_g\}$. For each generating unit, the minimum up and down limits $m^{\uparrow}_g$ and $m^{\downarrow}_g$ are randomly selected in such a way that $m^{\uparrow}_g-1$ and $m^{\downarrow}_g-1$ have Poisson distribution with parameter $4$. The initial status of generators at time period $t=0$ is found by solving a single period economic dispatch problem corresponding to the demand at time $t=1$. For each generating unit $g\in\mathcal{G}$, it is assumed that the initial status has been maintained exactly since time period $t=-t^{(0)}_g$, where $t^{(0)}_g$ has Poisson distribution with parameter $4$. For simplicity, all of the generating units with negative capacity are removed. Hourly load changes for the day-ahead at all buses are considered proportional to the numbers reported in \cite{khodaei2010transmission}. The changes in demand throughout the 24-hour planning horizon are reported in Table \ref{tab:UC_loadfactor}. For every time epoch, the corresponding demand factor at that time is multiplied by all loads in the system.

Table \ref{tab:UC_testcases} reports the results averaged over five Monte Carlo simulations for 24-hour scheduling. In this table, $k_{f}$ denotes the average round number of \cref{al:alg_1} at which the penalized relaxation produced a feasible solution with less than $10^{\smallMinus 6}$ per unit constraint violation. 

In order to evaluate the resulting feasible solutions from Algorithm \ref{al:alg_1} we solved an unpenalized semidefinite programming (SDP) relaxation of AC unit commitment by replacing the set $\Ccal$ in \cref{eq:UC_obj,eq:UC_set1,eq:UC_set2} with the cone of $|\Vcal|\times |\Vcal|$ Hermitian positive semidefinite matrices. The SDP relaxation offers a lower bound for the globally optimal cost of AC unit commitment, using which we can calculate the quality of our feasible solutions from Algorithm \ref{al:alg_1} through the formula 
\begin{equation}\small
\mathrm{GAP\%}=100\times\frac{\sum_{g,t}(c^{\mathrm{feasible}}_{g,t} -c^{\mathrm{SDP-lower-bound}}_{g,t})}{\sum_{g,t}c^{\mathrm{feasible}}_{g,t}},
\end{equation}
where $c^{\mathrm{feasible}}_{g,t}$ denotes the optimal cost value of the generating unit $g\in\Gcal$ at time $t\in\Tcal$ at round $50$ of the proposed sequential SOCP relaxation, and $c^{\mathrm{SDP-lower-bound}}_{g,t}$ denotes the cost values obtained from unpenalized SDP relaxation of \cref{eq:UC_obj,eq:UC_set1,eq:UC_set2}. 
The parameter $t(s)$ reports the average run time of all $50$ rounds of Algorithm \ref{al:alg_1} in seconds. The initial point of Algorithm \ref{al:alg_1} for all of the experiments is chosen as $\vbf_{\smbullet,t}^{\sm{0}{5}}\!=\!\onebf_{|\Vcal|}$, $\sbf_{\smbullet,t}^{\sm{0}{5}}\!=\!\pbf_{\mathrm{min}}$, $\vec{\sbf}_{\smbullet,t}^{\,\sm{0}{5}} = \diagrm\{\vec{\Cbf}\vbf_{\smbullet,t}^{\sm{0}{5}}{\vbf_{\smbullet,t}^{\sm{0}{5}}}^{\hspace{-0.0cm}\ast}\vec{\Ybf}^{\ast}\}$, 
$~\cev{\sbf}_{\smbullet,t}^{\,\sm{0}{5}} = \diagrm\{\cev{\Cbf}\vbf_{\smbullet,t}^{\sm{0}{5}}{\vbf_{\smbullet,t}^{\sm{0}{5}}}^{\hspace{-0.0cm}\ast}\cev{\Ybf}^{\ast}\}$, and $\xbf_{\smbullet,t}^{\,\sm{0}{5}}$ is set to the initial status of the generators, for all $t\in\Tcal$. 

For all of the random experiments, Algorithm \ref{al:alg_1} successfully finds a fully feasible operating point. Moreover, the reported gaps in Table \ref{tab:UC_testcases} demonstrate the effectiveness of our method in solving large instances of AC unit commitment. Changes in the resulting cost values with respect to the round numbers for one of the random experiments of each benchmark case are illustrated in \cref{plt:plot_feas_cost}.

\begin{table}
\vspace{0.22cm}
\caption{Hourly Demand Factor.}
\scriptsize
	\centering
    	 \begin{tabular}{|@{\;}c@{\;}|@{\;}c@{\;}|@{\;}c@{\;}|@{\;}c@{\;}|l}
\cline{1-4}
 Hour & Demand Factor & Hour & Demand Factor &  \\ \cline{1-4}

12:00 AM    & 0.6843 & 12:00 PM & 0.9460 &  \\ \cline{1-4}
01:00 AM    & 0.6451 & 01:00 PM & 0.9516 &  \\ \cline{1-4}
02:00 AM    & 0.6198 & 02:00 PM & 0.9721 &  \\ \cline{1-4}
03:00 AM    & 0.6044 & 03:00 PM & 0.9992 &  \\ \cline{1-4}
04:00 AM    & 0.6057 & 04:00 PM & 1.0000 &  \\ \cline{1-4}
05:00 AM    & 0.6269 & 05:00 PM & 0.9638 &  \\ \cline{1-4}
06:00 AM    & 0.6773 & 06:00 PM & 0.9608 &  \\ \cline{1-4}
07:00 AM    & 0.6937 & 07:00 PM & 0.9271 &  \\ \cline{1-4}
08:00 AM    & 0.7297 & 08:00 PM & 0.9270 &  \\ \cline{1-4}
09:00 AM    & 0.8084 & 09:00 PM & 0.9089 &  \\ \cline{1-4} 
10:00 AM    & 0.8930 & 10:00 PM & 0.7654 &  \\ \cline{1-4}
11:00 AM    & 0.9223 & 11:00 PM & 0.7641 &  \\ \cline{1-4}
	\end{tabular}
	\label{tab:UC_loadfactor}
\end{table}

         
	
\begin{table}
\caption{The performance of the proposed sequential penalized SOCP relaxation for 24-hour scheduling of IEEE benchmark systems.}
	\centering
	\scriptsize
\begin{tabular}{ |c||c|c|c|c|c|}
		\hline
		\multirow{2}{*}{Test Case} & \multicolumn{5}{c|}{SOCP}\\
		\cline{2-6}
 		& $\mu$ & $\alpha$ & $k_{f}$ & GAP(\%) & t(s) \\
		\hline
        \hline
		case57 & 1e0 & 1 & 1 & 0.00 & 603.0 \\
		\hline
		case118 & 1e0 & 10 & 1 & 2.27 & 1537.5 \\
        \hline
		case300 & 1e1 & 10 & 12.4 & 5.52 & 4010.0 \\
		\hline
	\end{tabular}
	\label{tab:UC_testcases}
\end{table}
\tikzset{
    master/.style={
        execute at end picture={
            \coordinate (lower right) at (current bounding box.south east);
            \coordinate (upper left) at (current bounding box.north west);
        }
    },
    slave/.style={
        execute at end picture={
            \pgfresetboundingbox
            \path (upper left) rectangle (lower right);
        }
    }
}
\pgfplotsset{every tick label/.append style={font=\footnotesize}}

\begin{figure*}
\vspace{0.25cm}
\captionsetup[subfigure]{margin={0.84cm,0cm}}
\hspace{-0.25cm}
\begin{subfigure}[normal]{0.33\linewidth}
\centering
\begin{tikzpicture}[master, scale=1]
\centering
\begin{axis}[
		width = 1\textwidth,
		height = 0.7\textwidth,
        xmode= normal,
        xmax=50,
        xmin=0,
		ymode= normal,
        xlabel = {\small{Rounds}},
        xlabel style={yshift=0.27cm},
        ylabel = {\small{Cost values}},
        ylabel style={yshift=-0.5cm},
        ytick = {1e6, 1.1e6, 1.2e6, 1.3e6},
        xtick = {0,10,20,30,40},
    	grid style = dashed,
        legend pos = north east,
		legend cell align = {left},
        legend image post style = {scale=0.7},
        legend style={font=\scriptsize, inner xsep=0.2pt, inner ysep=0.1pt},
]
\addplot[color=blue,line width=0.25mm] table[x=iter,y=objective,col sep=comma] {57.csv};
\end{axis}
\end{tikzpicture}
\vspace{-2mm}
\caption{}
\end{subfigure}
\begin{subfigure}[normal]{0.33\linewidth}
\centering
\begin{tikzpicture}[slave, scale=1]
\centering
\begin{axis}[
		width = 1\textwidth,
		height = 0.7\textwidth,
        xmode=normal,
		ymode=normal,
        xmax=50,
        xmin=0,
        xlabel = {\small{Rounds}},
        xlabel style={yshift=0.27cm},
        ylabel = {\small{Cost values}},
        ylabel style={yshift=-0.5cm},
        ytick = {2.4e6, 2.8e6, 3.2e6, 3.6e6},
        xtick = {0,10,20,30,40},
    	legend pos=south east,
    	grid style=dashed,
		legend cell align={left},
        legend image post style={scale=0.7},
        legend style={font=\scriptsize, inner xsep=0.2pt, inner ysep=0.1pt},
]
\addplot[color=blue,line width=0.25mm] table[x=iter,y=objective,col sep=comma] {118.csv};
\end{axis}
\end{tikzpicture}
\vspace{-2mm}
\caption{}
\end{subfigure}
\begin{subfigure}[normal]{0.33\linewidth}
\centering
\begin{tikzpicture}[slave, scale=1]
\centering
\begin{axis}[
		width = 1\textwidth,
		height = 0.7\textwidth,
        xmode=normal,
        xmax=50,
        xmin=0,
		ymode=normal,
        xlabel = {\small{Rounds}},
        xlabel style={yshift=0.27cm},
        ylabel = {\small{Cost values}},
        ylabel style={yshift=-.5cm},
        xtick = {0,10,20,30,40},
    	legend pos=south east,
    	grid style=dashed,
		legend cell align={left},
        legend image post style={scale=0.7},
        legend style={font=\scriptsize, inner xsep=0.2pt, inner ysep=0.1pt},
]
\addplot[color=blue,line width=0.25mm] table[x=iter,y=objective,col sep=comma] {300.csv};
\end{axis}
\end{tikzpicture}
\vspace{-2mm}
\caption{}
\end{subfigure}
\vspace{-2mm}
\caption{Convergence behavior of the proposed sequential penalized SOCP relaxation. The resulting cost values per round numbers are shown for one of the random experiments on each of the (a) IEEE 57 bus; (b) IEEE 118 bus; (c) IEEE 300 bus systems.}
\vspace{-0.3cm}
\label{plt:plot_feas_cost}
\end{figure*}
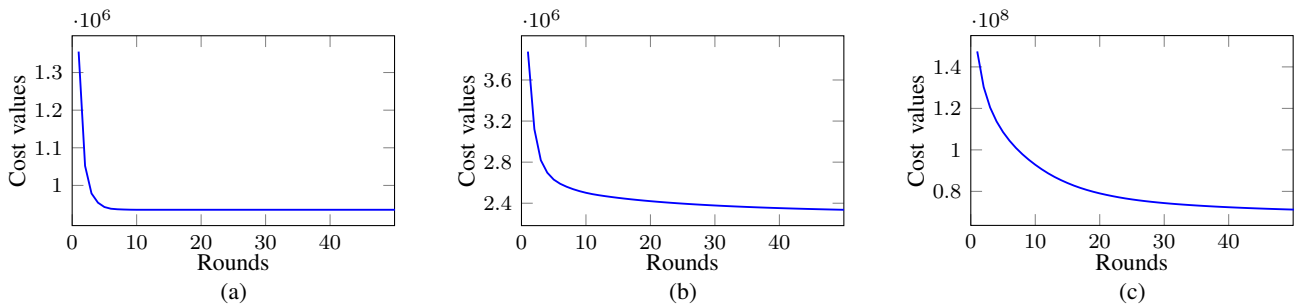

\section{Conclusions}
In this work, a sequential convex relaxation method is introduced for solving unit commitment with AC transmission constraints. We first, develop a second-order cone programming (SOCP) relaxation to convexity AC unit commitment problems. We then incorporate a penalty term into the objective of the proposed SOCP relaxation in order to find feasible solutions for the original non-convex AC unit commitment. 
The proposed penalized SOCP relaxations can be solved sequentially, to find feasible and near-globally optimal points.
The experimental results on IEEE 57 bus, IEEE 118 bus, IEEE 300 bus systems demonstrate the effectiveness of the proposed approach in solving challenging instances of AC unit commitment.

\bibliographystyle{IEEEtran}
\bibliography{IEEEabrv,egbib}
\end{document}